\documentclass[10pt]{amsart}

\usepackage{amssymb, latexsym, amsfonts, pigpen, enumitem, mathrsfs, amsthm, bbm, color}
\usepackage[matrix,arrow,curve]{xy}
\usepackage[plainpages,  urlcolor=blue]{hyperref}
\usepackage{graphicx}   

\newtheorem{theorem}[equation]{Theorem}
\newtheorem{lemma}[equation]{Lemma}
\newtheorem{corollary}[equation]{Corollary}

\theoremstyle{definition}
\newtheorem{definition}[equation]{Definition}
\newtheorem{example}[equation]{Example}
\newtheorem{examples}[equation]{Examples}

\theoremstyle{remark}
\newtheorem{remark}[equation]{Remark}
\numberwithin{equation}{section}

\makeatletter
\def\blfootnote{\gdef\@thefnmark{}\@footnotetext}
\makeatother

\newcommand{\GW}{{\mathrm{GW}}}
\newcommand{\SH}{{\mathrm{SH}}}
\newcommand{\Sm}{{\mathrm{Sm}}}

\newcommand{\ob}{{\mathrm{ob}}}
\newcommand{\W}{\mathrm{W}}

\newcommand{\MW}{\mathrm{MW}}

\newcommand{\A}{{\mathbb{A}}}
\newcommand{\C}{{\mathbb{C}}}
\newcommand{\R}{{\mathbb{R}}}
\newcommand{\F}{{\mathbb{F}}}
\newcommand{\Z}{{\mathbb{Z}}}
\newcommand{\HH}{{\mathbb{H}}}
\newcommand{\PP}{{\mathbb{P}}}

\newcommand{\QQ}{{\mathbb{Q}}}

\newcommand{\struct}{\mathcal{O}}

\newcommand{\GL}{{\mathrm{GL}}}
\newcommand{\SL}{{\mathrm{SL}}}
\newcommand{\GI}{\mathrm{GI}}
\newcommand{\Spec}{\operatorname{Spec}}
\newcommand{\tr}{\operatorname{tr}}
\newcommand{\Zar}{{\mathrm{Zar}}}

\newcommand{\Q}{\mathrm{Q}}
\newcommand{\chark}{\operatorname{char}}
\newcommand{\disc}{\delta}
\newcommand{\odeg}{\mathrm{deg_{GW}}}
\newcommand{\CHW}{{\widetilde{\mathrm{CH}}}}

\newcommand{\0}{{\emptyset}}
\newcommand{\sK}{{\mathcal{K}}}
\newcommand{\sO}{{\mathcal{O}}}

\newcommand{\del}{{\partial}}
\newcommand{\T}{{\mathfrak{T}}}

\begin{document}

\title{Combing a hedgehog over a field}

\author{Alexey Ananyevskiy}
\address{LMU M\"unchen \\ Mathematisches Institut \\ Theresienstr. 39 \\ 80333 M\"unchen \\ Germany}
\email{alseang@gmail.com}

\author{Marc Levine}
\address{
	Universit\"at Duisburg-Essen\\
	Fakult\"at Mathematik\\
	Thea-Leymann-Stra{\ss}e 9\\
	45127 Essen\\
	Germany}
\email{marc.levine@uni-due.de}

\date{July 1, 2024}

\begin{abstract} We investigate the question of the existence of a non-vanishing section of the tangent bundle on a smooth affine quadric hypersurface $Q^o$ over a given perfect field $k$. In case $Q^o$ admits a $k$-rational point, we give a number of necessary and sufficient conditions for such existence. We apply these conditions in a number of examples, including the case of the {\em algebraic $n$-sphere} over $k$, $S^n_k\subset \A^{n+1}_k$, defined by the equation $\sum_{i=1}^{n+1}x_i^2=1$. \end{abstract}

\maketitle

\tableofcontents

\section{Introduction}

\blfootnote{This paper is part of a project that has received funding from the European Research Council (ERC) under the European Union's Horizon 2020 research and innovation programme (grant agreement No. 832833).\\ 
\includegraphics[scale=0.08]{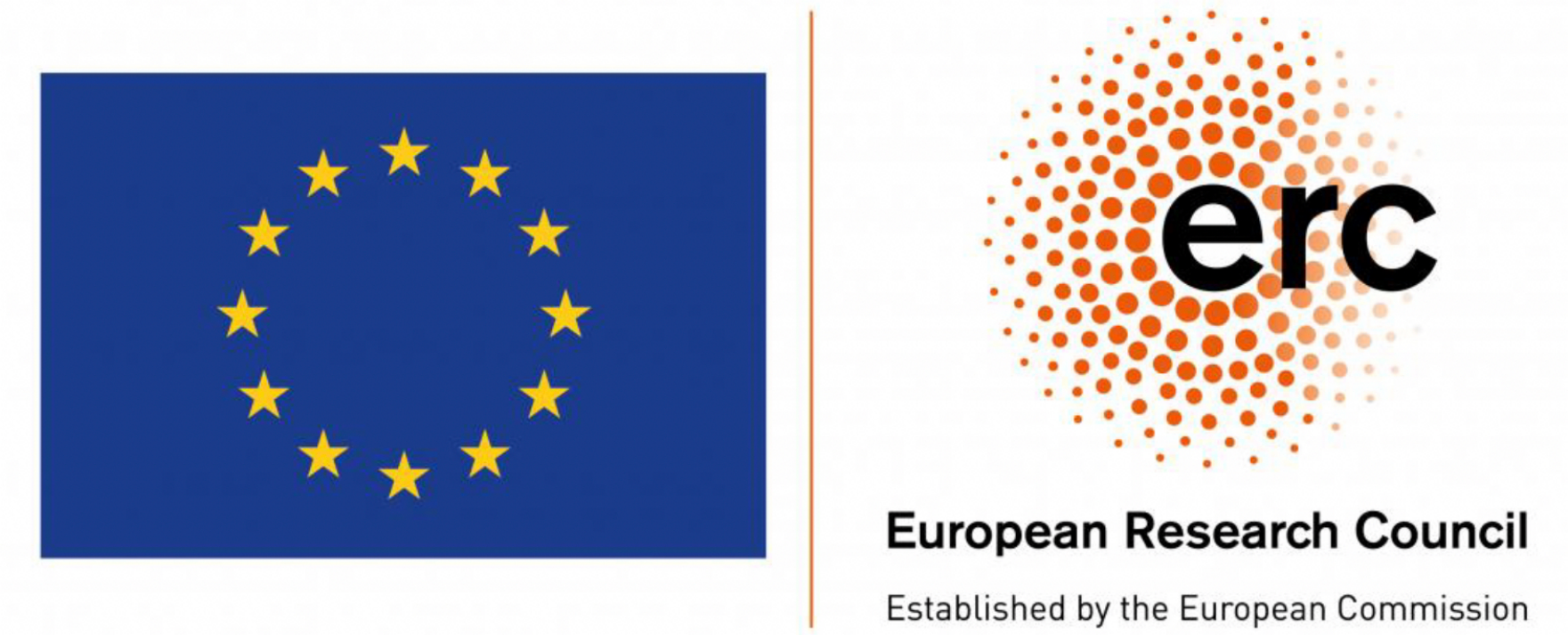}}

It is an elementary but nonetheless beautiful result found in nearly all introductory courses in differential topology, that for all $n\ge1$,  the tangent bundle $T_{S^{2n}}$ does not admit a non-vanishing section. One proof  uses the Gau{\ss}-Bonnet theorem to show that Euler class of $T_{S^{2n}}$ is non-zero by computing its degree as the Euler characteristic of $S^{2n}$, namely $2$, while the existence of a non-vanishing section would force the Euler class to vanish. For the odd dimensional case, the Euler characteristic vanishes, hence the  Euler class vanishes as well; one can also easily write down explicitly a non-vanishing section of $T_{S^{2n+1}}$. 

Writing the $n$-sphere $S^n$ as the hypersurface in $\R^{n+1}$ defined by the equation $\sum_{i=1}^{n+1}x_i^2=1$, one can ask the corresponding question in the algebro-geometric setting: let $k$ be a field of characteristic $\neq2$ and let $S^n_k\subseteq\A^{n+1}_k$ be the hypersurface defined by  the equation $\sum_{i=1}^{n+1}x_i^2=1$.  Does the tangent bundle $T_{S^n_k}$ admit a non-vanishing section?\footnote{To avoid any possible misunderstanding, for $E\to X$ a vector bundle on a $k$-variety $X$, a section $s:X\to E$ is said to be non-vanishing if the  scheme-theoretic intersection of $s(X)$ with the zero-section of $E$ is the empty scheme. Equivalently, letting $\bar{k}$ denote the algebraic closure of $k$, the set of $\bar{k}$ points $x$ of $X$ with $s(x)=0$ is empty.}

This question for $S^2_{\QQ_p}$ was originally raised by Umberto Zannier (see Remark~\ref{rem:Unimodular} below for the original formulation by Zannier). He showed that $S^2_{\QQ_p}$ admits a non-vanishing vector field for odd $p$ and he asked if there is a non-vanishing vector field on the 2-sphere over $\QQ_2$, motivating our interest in the question of the existence of non-vanishing vector fields on $S^n_k$ for arbitrary $n$ and $k$. 

We give an essentially complete answer to this question;  in case $k$ is perfect, this is in fact a special case of the more general Theorem~\ref{IntroThm:Main} about smooth affine quadric hypersurfaces.

\begin{theorem}[see Theorem~{\hyperref[thm:hedgehog]{\ref*{thm:hedgehog}.(1,3)}} and Remark~\ref{rem:Explicit}] Let $k$ be a field of characteristic $\neq2$. 
	\begin{enumerate}
		\item If $n$ is odd, then $T_{S^n_k}$ admits a non-vanishing section.
		\item If $n>0$ is even, then $T_{S^n_k}$ admits a non-vanishing section if and only if $-1$ is in the subgroup of $k^\times$ generated by the non-zero values of the function $\sum_{i=1}^{n+1}x_i^2$ on $k^{n+1}$.
	\end{enumerate}
\end{theorem}

As the condition in (2) is not very explicit, we reformulate this as follows

\begin{corollary}[Corollary~\ref{cor:hedgehog}]\label{IntroThm2} Let $k$ be a field of characteristic $\neq2$. For $n>0$ even, 
	$T_{S^n_k}$ admits a non-vanishing section if and only if the equation
	\[
	\sum_{i=1}^{2n+1}x_i^2=-1
	\]
	has a solution with the $x_i\in k$.
\end{corollary}

\begin{examples} \label{IntoExamples} Let $S^n_k$ be as above.
	\begin{enumerate}
		\item Suppose that $\chark k =p>2$. Then $T_{S^n_k}$ has a non-vanishing section for all $n>0$.
		\item Suppose $k$ contains a $p$-adic field $\QQ_p$. Then  $T_{S^n_k}$ has a non-vanishing section for all $n>0$.
		\item Suppose that $k$ is a number field, and take $n>0$ to be even. Then $T_{S^n_k}$ has a non-vanishing section if and only if $k$ has no real embeddings.
	\end{enumerate}

	To see this we apply Corollary~\ref{IntroThm2}.	For (1), since $ \F_p\subseteq k$, it suffices to take $k=\F_p$. Since every element $x\in \F_p^\times$ is a sum of two squares \cite[Proposition~II.3.4]{Lam05}, the condition of Corollary~\ref{IntroThm2} is satisfied for all $n\ge1$. See also Remark~\ref{rem:Explicit} for an explicit non-vanishing section.
	
	For (2), we reduce as above to the case $k=\QQ_p$. If $p$ is odd, then by Hensel's lemma, each solution to 
	$ \sum_{i=1}^{2n+1}x_i^2=-1$ in $\F_p$ lifts to a solution in $\Z_p$, so the criterion is satisfied. For $p=2$, the class of a unit $u$ in $\Z_2$ modulo squares is given by the image of $u$ in $(\Z/8)^\times$, so it suffices to write $7$ as the sum of $\le 5$ squares in $\Z_2$ and it turns out that $4$ squares are enough: $7=1+1+1+4$. For those more intrinsically minded, one has the general result that {\em every} non-degenerate quadratic form $\phi$ in at least 5 variables over a local field  has a non-trivial zero \cite[Theorem~VI.2.12]{Lam05}, which we apply to $\phi=\sum_{i=1}^5x_i^2$.
	
	For (3), it is clear that the equation $ \sum_{i=1}^{2n+1}x_i^2=-1$ has no solution in $k^{2n+1}$ if $k$ admits a real embedding. Conversely, we may use the Hasse-Minkowski principle for quadratic forms (see e.g., \cite[Hasse-Minkowski Principle~VI.3.1]{Lam05}) to see that $\sum_{i=1}^4x_i^2=-1$ has a solution in $k$ if $k$ is a purely imaginary number field. Indeed, it suffices to show that $\sum_{i=1}^4x_i^2=-1$ has a solution in $k_v$ for every place $v$ of $k$. This is clear if $v$ is an infinite place, as $k_v=\C$ by assumption. If $v$ is a finite place, then $k_v\supset\QQ_p$ for some prime $p$, and we have just seen that 
	$\sum_{i=1}^4x_i^2=-1$ has a solution in $\QQ_p$ for every prime $p$. 
\end{examples}

One can also ask about a general smooth affine quadric $Q^o\subseteq\A^{n+1}_k$, with $k$ a field of characteristic $\neq2$. Since every quadratic form over $k$ can be diagonalized, we may assume that  $Q^o$ is  defined by an equation of the form $q=1$, where $q=\sum_{i=1}^{n+1}a_ix_i^2\in k[x_1,\ldots, x_{n+1}]$, with $\prod_ia_i\neq0$. Here one has a result of essentially the same form as for $S^n_k$, with the extra condition that for even $n$, $Q^o(k)$ should be non-empty, that is, $q=1$ has a solution in $k$.

Let $D(q)$ be the set of non-zero values of $q$ on $k^{n+1}$, let $D(q)^2=\{a\cdot b\mid a,b\in D(q)\}\subseteq k^\times$,   and let $[D(q)]$, $[D(q)^2]$ be the subgroups of $k^\times$ generated by $D(q)$,  $D(q)^2$, respectively.

\begin{theorem}[Theorem~{\hyperref[thm:hedgehog]{\ref*{thm:hedgehog}.(1,3)}}]\label{IntroThm:Main} Let $k$ be a perfect field of characteristic $\neq 2$, let $q=\sum_{i=1}^{n+1}a_ix_i^2$ with $a_1,\ldots, a_{n+1}\in k^\times$ and let $Q^o\subseteq\A^{n+1}_k$ be the affine quadric hypersurface $q=1$.  
	\begin{enumerate}
		\item If $n$ is odd, then $T_{Q^o}$ has a non-vanishing section.
		\item Suppose $Q^o(k)\neq\0$. If $n>0$ is even, then $T_{Q^o}$ has a non-vanishing section if and only if $-1\in [D(q)]$.
	\end{enumerate}
\end{theorem}
\noindent 
If $Q^o(k)=\0$ and $n$ is even, we only have a necessary condition for the existence of a 
non-vanishing section of $T_{Q^o}$. 

\begin{theorem}[Theorem~{\hyperref[thm:hedgehog]{\ref*{thm:hedgehog}.(2)}}]\label{IntroThmNecessary} Let $k$, $q$ and $Q^o$ be as above. If $n$ is even and $T_{Q^o}$ has a non-vanishing section, then $-\prod_{i=1}^{n+1}a_i\in [D(q)^2]$.
\end{theorem}
\noindent 
Since $a_i\in D(q)$ for each $i$, the above condition is the same as asking for $-a_i$ to be in $[D(q)^2]$ for {\em some} $i$. Note that  $Q^o(k)\neq\0$ if and only if $1\in D(q)$, so in case  $Q^o(k)\neq\0$, we have $[D(q)^2]=[D(q)]$, and $-1\in[D(q)]\Leftrightarrow -\prod_{i=1}^{n+1}a_i\in [D(q)]$.

Here is a version of Example~\ref{IntoExamples} for general $q$.

\begin{corollary}[Corollaries~\ref{cor:CohDim2},~\ref{cor:NumberField}, and~\ref{cor:LocalFieldEtc}] Let $k$, $q$ and $Q^o$ be as in Theorem~\ref{IntroThm:Main}.
	\begin{enumerate}
		\item Let $k=\F_{p^m}$  with $p>2$. Then $T_{Q^o}$ has a non-vanishing section for all $n>0$.
		\item Suppose $k$ is a non-Archimedean local field of characteristic zero,  the perfection of a local field of characteristic $p>2$, or the perfection of a function field of a curve over a finite field of characteristic $p>2$.  Then for $n$ odd, or $n\ge4$ even,  $T_{Q^o}$ has a non-vanishing section. If $n=2$, then $T_{Q^o}$ has a non-vanishing section if $Q^o(k)\neq\0$. 
		\item Suppose $k$ is a number field, $Q^o(k)\neq\0$ and $n>0$ is even. Then $T_{Q^o}$ has a non-vanishing section if and only if the equation $q=0$ has a non-trivial solution in $k_v$ for every real place $v$ of $k$. Equivalently, for each real embedding $\sigma:k\hookrightarrow\R$,  $\sigma(a_i)<0$  for some $i$. $T_{Q^o}$ also has a non-vanishing section if $n$ is odd. 
		\item Let $k$ be a perfect field of cohomological dimension $\le 2$. Suppose that $n$ is odd, or that $n>0$ is even and $Q^o(k)\neq\0$. Then $T_{Q^o}$ has a non-vanishing section. 
	\end{enumerate}
\end{corollary}

\begin{remark}[Uni-modular rows and uni-modular matrices]\label{rem:Unimodular} Let $q=\sum_{i=1}^{n+1}a_ix_i^2\in k[x_1,\ldots, x_{n+1}]$ be a quadratic form, defining $Q^o\subseteq \A^{n+1}_k$ as $V(q-1)$, and let  $R$ be the coordinate ring 
	\[
	R:=k[x_1,\ldots, x_{n+1}/(q-1).
	\]
	We are assuming that $q$ is non-degenerate, that is,  $\prod_ia_i\neq0$, and that $n\ge1$.
	
	Let $\nabla(q)$ denote the gradient
	\[
	\nabla(q):=(\del q/\del x_1,\ldots, \del q/\del x_{n+1})
	\]
	and let 
	\[
	\tilde{\nabla}(q):=(a_1x_1,\ldots, a_{n+1}x_{n+1}),
	\]
	so $2\tilde{\nabla}(q)=\nabla(q)$. 
	
	We first assume that $2$ is invertible in $k$, so we can rewrite the 
	tangent-normal sequence for $Q^o\subseteq \A^{n+1}_k$ using $\tilde{\nabla}(q)$, as
	\begin{equation}\label{eqn:TangentNormal}
		0\to T_{Q^o}\to \sO_{Q^o}^{n+1}\xrightarrow{\tilde{\nabla}(q)^t} \sO_{Q^o}\to 0.
	\end{equation}
	where $\tilde{\nabla}(q)^t\in M_{n+1\times 1}(R)$ is the transpose of $\tilde{\nabla}(q)$.
	Since $Q^o$ is affine, we can rephrase everything in terms of $R$-modules, giving the exact sequence
	\[
	0\to \T_{Q^o}\to R^{n+1}\xrightarrow{\tilde{\nabla}(q)^t} R\to 0,
	\]
	with $\T_{Q^o}$ the $R$-module of global sections of $T_{Q^o}$. Since
	\[
	(x_1,\ldots,x_{n+1})\cdot \tilde{\nabla}(q)=1
	\]
	we may split the surjection by $(x_1,\ldots,x_{n+1}):R\to R^{n+1}$, exhibiting $\T_{Q^o}$ as a stably free $R$-module, and showing that $(x_1,\ldots,x_{n+1})$ is a {\em uni-modular row}, i.e.,   $(x_1,\ldots,x_{n+1})R$ is the unit ideal.
	
	It is straightforward to see that the stably free $R$-module $\T_{Q^o}$ is free if and only if there is a matrix $M\in \GL_{n+1}(R)$ with the first row $(x_1,\ldots,x_{n+1})$; by dividing the last row of $M$ by $\det M$, we may in fact take $M$ to have $\det M=1$, so  $\T_{Q^o}$ is a  free $R$-module if and only if there is a uni-modular matrix $M$ over $R$ with the first row $(x_1,\ldots,x_{n+1})$.
	
	More generally, we may take $k$ to be an arbitrary commutative ring (even with $2$ not a unit), and let $(a_{ij})_{1\le i,j\le n+1}\in \GL_{n+1}(k)$ be an invertible symmetric matrix. Let
	\begin{gather*}
		q:=\sum_{i,j=1}^{n+1}a_{ij}x_ix_j, \ R:=k[x_1,\ldots, x_{n+1}]/(q-1),\ Q^o:=\Spec R,\\
		\tilde{\nabla}(q):=(\sum_ja_{1j}x_j,\sum_ja_{2j}x_j,\ldots, \sum_ja_{{n+1}j}x_j).
	\end{gather*}
	Then the map $R^{n+1}\xrightarrow{\tilde{\nabla}(q)^t} R$ is surjective, and  we may define a stably free $R$-module $\T_{Q^o/k}$ by the exact sequence 
	\begin{equation}\label{eqn:DefTQ}
		0\to \T_{Q^o/k}\to R^{n+1}\xrightarrow{\tilde{\nabla}(q)^t} R\to 0.
	\end{equation}
	Since  $(x_1,\ldots,x_{n+1})\cdot \tilde{\nabla}(q)=1\in R$, $(x_1,\ldots,x_{n+1})$ is a uni-modular row over $R$, and the $R$-module $\T_{Q^o/k}$ is free if and only if there is a matrix $M\in \SL_{n+1}(R)$ with first row $(x_1,\ldots,x_{n+1})$. Furthermore, if $2\in k^\times$, then $Q^o$ is smooth over $k$, and $\T_{Q^o/k}$ is the $R$-module of global sections of the relative tangent bundle $T_{Q^o/k}\to Q^o$ of $Q^o$ over $\Spec k$. 
	
	Note that, for $n=2$ and $2\in k^\times$, the existence of a uni-modular matrix over $R$ with the first row $(x_1,x_2,x_3)$ is equivalent to the existence of a non-vanishing section of $T_{Q^o/k}$. Indeed, $T_{Q^o/k}$ admits a non-vanishing section if and only if we can write $\T_{Q^o/k}\cong R\oplus P$, with $P$ a rank one projective $R$-module, which yields the isomorphism $\bigwedge^2_R\T_{Q^o/k}\cong P$. The exact sequence \eqref{eqn:DefTQ} gives an isomorphism $\bigwedge^2_R\T_{Q^o/k}\cong R$, so $P\cong R$ and thus $\T_{Q^o/k}\cong R^2$. This gives the existence of $M\in \SL_3(R)$ with first row $(x_1,x_2,x_3)$. 
	
	The original form of Zannier's question was in the following terms: taking $q:=x^2+y^2+z^2\in k[x,y,z]$, for which fields $k$ is the uni-modular row $(x,y,z)$ over $R$ the first row in a uni-modular matrix over $R$? He showed this was the case for $k=\QQ_p$, $p$ odd, and asked about $k=\QQ_2$.  Note that one can just as well ask the question for $k$ an arbitrary commutative ring and general $q$ as above; our results only give a criterion for a positive answer to this question in case $k$ is a perfect field with $2$ invertible and $n=2$. Noting our positive answer for $k=\QQ_2$, Zannier asked in a recent private communication about the case $k=\Z_2$.  \end{remark}

\begin{remark}[Explicit sections]\label{rem:Explicit} 1. In case $n$ is odd or if the quadratic form $q:=\sum_{i=1}^{n+1}a_ix_i^2$ is isotropic over $k$ (i.e., $q=0$ has a non-trivial solution in $k$), then one can write down explicit non-vanishing sections of $T_{Q^o}$. 
	
	For $n$ odd,  
	the tangent-normal sequence for $Q^o\subseteq \A^{n+1}_k$,
	\[
	0\to T_{Q^o}\to \sO_{Q^o}^{n+1}\xrightarrow{(2a_1x_1,\ldots, 2a_{n+1}x_{n+1})^t}\sO_{Q_o}\to 0,
	\]
	says a section $s$ of $T_{Q^o}$ is given by an $n+1$-tuple of regular functions $(s_1,s_2,\ldots, s_{n+1})$ with $\sum_{i=1}^{n+1}a_ix_is_i=0$. One can take 
	\[
	s=(a_2x_2,-a_1x_1,\ldots, a_{n+1}x_{n+1}, -a_nx_n)
	\]
	which is clearly non-vanishing. 
	
	This is a special case of the following general result. Let $A$ be   a commutative ring,  let  $(a_1,\ldots, a_{2m})$ be a uni-modular row  in $A^{2m}$ (i.e.,  $a_1,\ldots, a_{2m}$ generate the unit ideal in $A$), and let  $M$ be the stably free $A$-module defined by the exactness of the sequence
	\[
	0\to M\to A^{2m}\xrightarrow{(a_1,\ldots, a_{2m})^t} A\to0.
	\]
	Then $M$ admits the free rank one  summand defined by
	\[
	0\to A\xrightarrow{(-a_2,a_1,\ldots, -a_{2m}, a_{2m-1})}M\subset A^{2m}.
	\]
	
	Now take $n$ even and suppose $q$ is isotropic. Then after a $k$-linear change of coordinates and change of notation, we may assume that 
	\[
	q=2x_1x_2+\sum_{i=3}^{n+1}a_ix_i^2
	\]
	(see e.g. \cite[Theorem~I.3.4]{Lam05}). In this case, the tangent-normal sequence for $Q^o\subseteq \A^{n+1}_k$ is
	\[
	0\to T_{Q^o}\to \sO_{Q^o}^{n+1}\xrightarrow{(2x_2, 2x_1,2a_3x_3,\ldots, 2a_{n+1}x_{n+1})^t} \sO_{Q_o}\to 0.
	\]
	Letting
	\[
	s=(0, a_3x_3, -x_1, a_5x_5, -a_4x_4,\ldots, a_{n+1}x_{n+1}, -a_nx_n)
	\]
	gives a section of $T_{Q^o}$ with $s=0$ given by $Q^o\cap (x_1=x_3=\ldots=x_{n+1}=0)$, which is clearly empty. 
	
	In particular, let $k$ be a field of characteristic $p>2$. In a finite field $-1$ is a sum of two squares \cite[Proposition~II.3.4]{Lam05} whence the quadratic form $x_1^2+x_2^2+\ldots+x_{n+1}^2$ is isotropic over $k$ provided that $n\ge 2$. Hence the tangent bundle $T_{S^n_k}$ to the algebraic $n$-sphere over $k$ admits a non-vanishing section for every $n\ge 1$.\\[5pt]
	2. Our main results for  even $n$ and $q$ anisotropic only give criteria for the  existence of a non-vanishing section, without giving an explicit formula. In the case that Zannier had asked about originally, $S^2_{\QQ_2}$, Peter M\"uller \cite{M24}, noting our existence result and following a suggestion of Zannier,   found an explicit trivialization of $T_{S^2_{\QQ_2}}$. We quote from his private communication: \\[2pt]
	``Indeed, some sophisticated computations eventually gave an explicit
	example over $\QQ_2$ (in fact even over $\QQ(\sqrt{-7})$), \ldots'' \\[2pt]
	Here is M\"uller's example.  Let  $R=\QQ_2[x,y,z]/(x^2+y^2+z^2-1)$,  the coordinate ring of $S^2_{\QQ_2}$. The polynomial $T^2-T+2$ has two roots in $\Z_2$, exactly one of which, which we denote by  $\omega$,   is a unit in $\Z_2$. In particular,  $2-\omega$ is also a unit in $\Z_2$. M\"uller gives his example in the form of a $3\times 3$ matrix over $R$ with determinant $2-\omega$ and first row  $(x,y,z)$. The explicit matrix is
	\[
	\begin{pmatrix}
		x&y&z\\
		-y+z+1&(1-\omega)x+y+(1+\omega)z+\omega&
		-x-y+2z+(1-\omega)\\
		\omega y +(2-\omega)z&(1-2\omega)x+(1+\omega)y+3z+1&-2x+(2-\omega)z-\omega
	\end{pmatrix}\lower13pt\hbox{.}
	\]
	Let $\lambda_i$ be the dot product of $(x,y,z)$ with the $i$th row. Noting that
	$(x,y,z)\cdot (x,y,z)=1$ 
	in $R$,  this gives the following two independent non-vanishing sections of $T_{S^2_{\Q_2}}$:
	\begin{gather*}
	s_1(x,y,z)=(-y+z+1,(1-\omega)x+y+(1+\omega)z+\omega,
	-x-y+2z+(1-\omega))-\lambda_2\cdot(x,y,z),\\
	s_2(x,y,z)=(\omega y +(2-\omega)z,(1-2\omega)x+(1+\omega)y+3z+1,-2x+(2-\omega)z-\omega)-\lambda_3\cdot(x,y,z).
	\end{gather*}
	M\"uller notes that this also works over $\QQ(\sqrt{-7})$, where we take $\omega$ to be either of the two roots of $T^2-T+2$ in $\QQ(\sqrt{-7})$.
	
	In addition, M\"uller's example gives a positive answer to Zannier's question over $\Z_2$ instead of $\QQ_2$, just divide the last row by the determinant $2-\omega\in \Z_2^\times$. 
\end{remark}

\begin{remark}[Some non-examples] Suppose $n$ is even. We have already seen in Remark~\ref{rem:Explicit} that $T_{Q^o}$ has a non-vanishing section if $q$ is isotropic over $k$. On the other hand, $q$ being isotropic over $k$ implies that $Q^o$ has a $k$-rational point \cite[Theorem~I.3.4.(3)]{Lam05}, so if $Q^o(k)=\0$, then $q$ is anisotropic over $k$ and we do not have any explicit method for finding a (possible) non-vanishing section of $T_{Q^o}$. Moreover, Theorem~\ref{IntroThmNecessary} is our only result that considers the case $n$ even and $Q^o(k)=\0$,  and it only gives us  a necessary condition for $T_{Q^o}$ to have a non-vanishing section. Here is a series of examples that are not covered by any of our results.
	
	Take $k=\QQ_p$ with $p>2$. Let $u\in \Z_p^\times$ be a non-square modulo $p$, and let  $q=ux_1^2+px_2^2-upx_3^2$. It follows from \cite[Theorem~VI.2.2]{Lam05}  that $q-x_0^2$ is anisotropic over $\QQ_p$, hence $Q^o(\QQ_p)=\0$ and also $q$ is anisotropic over $\QQ_p$. Moreover $-(u\cdot p\cdot(-up))=1$ in $\QQ_p^\times/\QQ_p^{\times 2}$, so $-(u\cdot p\cdot(-up))$ is in $[D(q)^2]$, hence the necessary condition in Theorem~\ref{IntroThmNecessary} is satisfied.  We do not know whether $T_{Q^o}$ has a non-vanishing section in any of these cases.
	
	One final non-example. Take $k=\R$, $q=\sum_{i=1}^{n+1}-x_i^2$, with $n$ even, and let $R$ be the coordinate ring of $Q^o$. Then we have $Q^o(\R)=\0$.  By a theorem of Ojanguren-Parimala-Sridharan \cite[Theorem 3.2]{OPS86}, there is an $M\in \SL_{n+1}(R)$ with
	\[
	(x_1,\ldots, x_{n+1})=(1,0,\ldots, 0)M
	\]
	in other words, $(x_1,\ldots, x_{n+1})$ is the first row of the uni-modular matrix $M$. Thus the $R$-module $\T_{Q^o}$ is free, so $T_{Q^o}$ is a trivial vector bundle over $Q^o$, hence admits a non-vanishing section. Our results only yield the necessary condition 
	\[
	-(-1)^{n+1}=1\in|D(q)^2|, 
	\]
	which (fortunately!) is true in this case.  
\end{remark}

The main idea behind the proofs of our results is quite close to the case of the real spheres. We have already disposed of the case of odd $n$ in  Remark~\ref{rem:Explicit}. For $n>0$ even, we replace the Euler class $e_{top}(T_{S^n})\in H^n(S^n,\Z)$ with the Euler class $e(T_{Q^o})$ in the {\em Chow-Witt group} $\CHW^n(Q^o)$. For a smooth $k$-variety $X$ and a rank $r$ vector bundle $E$ on $X$, one has an Euler class $e(E)$ in the (twisted) Chow-Witt group $\CHW^n(X,\det^{-1}(E))$. In our case, the tangent-normal sequence for $T_{Q^o}$ gives a canonical isomorphism $\det T_{Q^o}\cong \sO_{Q^o}$, which induces an isomorphism $\CHW^n(Q^o, \det^{-1}T_{Q^o})\cong \CHW^n(Q^o)$ with the untwisted version of the Chow-Witt group, giving us our Euler class $e(T_{Q^o})\in\CHW^n(Q^o)$.  A fundamental result of Morel \cite[\S~8.2]{Mor12} says that for a smooth affine $k$-scheme $X$ of dimension $n$ over a perfect field $k$ and a rank $n$ vector bundle $E$ over $X$, $E$ admits a non-vanishing section if and only if the Euler class $e(E)$ vanishes (in this form, the result also relies on work of Asok-Fasel \cite{AF16} and Asok-Hoyois-Wendt \cite{AHW17}, see Theorem~\ref{thm:Morel_criterion} for the discussion).

Since $Q^o$ is not proper over $k$, we do not have a nice analog of the Gau{\ss}-Bonnet theorem for 
$Q^o$, so we pass to its projective closure $Q\subseteq\PP^{n+1}$, defined by the equation $\sum_{i=1}^{n+1}a_ix_i^2=x_0^2$, and let $Q^\infty\subseteq Q$ be the hyperplane section defined by $x_0=0$. 

Let $\GW(k)$ denote the Grothendieck-Witt ring of (virtual) regular quadratic forms over $k$. For $p:X\to \Spec k$ a smooth projective variety, we have the pushforward map
\[
p_*:\CHW_0(X)\to \CHW_0(\Spec k)=\GW(k),
\]
which we denote by $\odeg:\CHW_0(X)\to\GW(k)$; we call this the {\em quadratic degree map}. $X$ has a {\em quadratic Euler characteristic} $\chi(X/k)\in \GW(k)$ and we have a quadratic Gau{\ss}-Bonnet theorem (\cite[Theorem 4.6.1]{DJK21}, \cite[Theorem 5.3]{LR20}): letting $T_X$ be the tangent bundle of $X$, we have
\[
\odeg(e(T_X))=\chi(X/k),
\]
so we are all set up to argue as in differential topology.

Getting back to our quadrics, let us first assume that $Q^o$ has a $k$-rational point. We show in Section~2 that  $e(T_{Q^o})=0$ if and only if $\chi(Q/k)$ is in the subgroup $\odeg(\CHW_0(Q^\infty))\subseteq \GW(k)$, and we have an explicit expression for $\chi(Q/k)$:
\[
\chi(Q/k)=\langle 2, 2\prod_ia_i\rangle+\frac{n}{2}\langle1,-1\rangle,
\]
where $\langle a,b\rangle$ is the quadratic form $ax^2+by^2$, and $m\cdot \langle a,b\rangle$ is the quadratic form $\sum_{i=1}^max_i^2+by_i^2$.

Putting this together, we see that $T_{Q^o}$ admits a non-vanishing section if and only if 
$\langle 2, 2\prod_ia_i\rangle+\frac{n}{2}\langle1,-1\rangle$ is in $\odeg(\CHW_0(Q^\infty))\subseteq \GW(k)$.  We conclude Section~2 by using this criterion to handle the cases discussed in Examples~\hyperref[IntoExamples]{\ref*{IntoExamples}.(1,2)} above.   

The next step is to use the theory of quadratic forms to rephrase the condition 
\[
\langle 2, 2\prod_ia_i\rangle+\frac{n}{2}\langle1,-1\rangle\in \odeg(\CHW_0(Q^\infty))\subseteq\GW(k)
\]
in terms of the subgroups $[D(q)]$ and $[D(q)^2]$. This  is done in Section~3, relying on properties of Scharlau's transfer, Knebusch's norm principle and basic facts about Pfister forms and Pfister neighbors. We apply these tools in Section~4 to give our main results, which yield criteria that are much easier to apply than the one derived in Section~2. We conclude by using this to compute the remaining examples described above; the case in which $Q^o$ does not have a $k$-rational point is trickier to handle, and we are only able to arrive at the necessary condition stated in Theorem~\ref{IntroThmNecessary}.

\textbf{Acknowledgements.} The work of the first author is supported by the DFG research grant AN 1545/4-1 and DFG Heisenberg grant AN 1545/1-1. Alexey Ananyevskiy would like to thank Fabien Morel who drew his attention to the problem of the existence of non-vanishing vector fields on a 2-sphere over the dyadic numbers. Fabien Morel in turn learned about this problem from Umberto Zannier some time ago. This paper grew out of this particular question during the first author's research visit to the University of Duisburg-Essen.

The second author is supported by the ERC grant QUADAG. Marc Levine would like to thank Christopher Deninger, who passed on the problem of the existence of non-vanishing vector fields on a 2-sphere over the dyadic numbers from Umberto Zannier in 2020, as well as Umberto Zannier himself for subsequent discussions. We would like to thank the referee for their very helpful comments and suggestions, and we thank Peter M\"uller for kindly allowing us to present his example here.

Throughout the paper we employ the following  notations.
\\[5pt]
	\begin{tabular}{l|l}
	$k$ & a perfect field with $\chark k\neq 2$ \\
	$\Sm_k$& the category of smooth, separated, finite-type $k$-schemes\\
	$T_X$ & the tangent bundle of $X\in \Sm_k$\\
	$X(F)$ & the set of rational points of $X_F$ for $X\in\Sm_k$ and a field extension $F/k$\\
	$\GW(F)$ & the Grothendieck-Witt ring of (virtual) regular quadratic forms over a field $F$ \\
	$\W(F)$ & the Witt ring, the quotient of $\GW(F)$ by the hyperbolic forms \\
	$\GI(F)$ & the ideal in $\GW(F)$ generated by the even-dimensional forms\\
	$I(F)$ & the image in $\W(F)$ of $\GI(F)$\\	
	$F^\times$ & the multiplicative group of non-zero elements of the field $F$\\
	$\langle a_1,a_2,\dots,a_n\rangle$ & the quadratic form $a_1x_1^2+a_2x_2^2+\dots + a_nx_n^2$ \\
	\end{tabular}

\section{Recollections on Chow-Witt groups and a computational criterion}

\begin{definition}
	Assume $k$ to be a perfect field. We will use the Chow-Witt groups, also known as Chow groups of oriented cycles, that were introduced in \cite{BM00}. These groups provide refined cohomological obstructions to the existence of non-vanishing sections of algebraic vector bundles. We refer the reader to the expositions in \cite{Fas20, Deg23} and \cite[Sections~2,3]{AF16} for the properties of these groups that we list below.
	
	We recall from \cite[Chapter~2]{Mor12} the {\em Milnor-Witt $K$-theory sheaves} $\sK_n^{\MW}$, $n\in\Z$. These are Nisnevich sheaves of abelian groups on $\Sm_k$, with products $\sK_n^{\MW}\times \sK_m^{\MW}\to \sK_{n+m}^{\MW}$ making the graded object $\sK_*^{\MW}:=\oplus_{n\in\Z}\sK_n^{\MW}$ into a sheaf of associative, unital, graded rings on $\Sm_k$. Given $X\in \Sm_k$ and a line bundle $L$ on $X$, we have the $L$-twisted version $\sK_n^{\MW}(L)$, giving a Nisnevich sheaf on $\Sm_k/X$.
	Letting $\mathcal{GW}$ denote the Nisnevich sheaf  of Grothendieck-Witt rings, there is a canonical isomorphism  $\sK_0^{\MW}\cong \mathcal{GW}$. 
	For a  field $F$, and $L$ a dimension one $F$-vector space, we write  $K_n^{\MW}(L)(F)$ for $\sK_n^{\MW}(L)(\Spec F)$. 

	For a smooth variety $X$ over $k$, a line bundle $L$ over $X$ and an integer $n\ge0$, the \textit{Chow-Witt group} $\CHW^n(X,L)$ is defined as:
	\[
	\CHW^n(X,L) =H_\Zar^n(X; \sK^{\MW}_n(L)).
	\]
	We will also use the homological notation with
	\[
	\CHW_n(X,L) = \CHW^{d-n}(X,L\otimes \omega_X),
	\]
	where $d=\dim X$ and $\omega_X$ is the canonical bundle of $X$. We put
	\[
	\CHW^n(X):=\CHW^n(X,\struct_X),\quad \CHW_n(X):=\CHW_n(X,\struct_X).
	\]
	Chow-Witt groups have the following properties that we will use below.
	\begin{enumerate}
		\item $\CHW^n(X,L)$ is canonically identified by \cite[Theorem~5.47]{Mor12} with the $n$-th cohomology group of the Rost-Schmid complex 
		\begin{multline*}
			\bigoplus_{x\in X^{(0)}} K^{\MW}_n(L_x\otimes \omega_{x/X})(k(x)) \to \bigoplus_{x\in X^{(1)}}  K^{\MW}_{n-1}(L_x\otimes \omega_{x/X})(k(x))\to \\ \to \dots \to 
			\bigoplus_{x\in X^{(d)}}  K^{\MW}_{n-d}(L_x\otimes \omega_{x/X})(k(x)),
		\end{multline*}
		where the sums are taken over all the points of $X$ of the respective codimension, $L_x$ is the restriction of $L$ to $x$, $\omega_{x/X}$ is the determinant of the normal bundle for the embedding $x\to X$ and $d=\dim X$.
		
		\item
		For a line bundle $L'$ over $X$ there is a canonical isomorphism \cite[Remark~5.13]{Mor12}
		\[
		\CHW^n(X,L\otimes (L')^{\otimes 2}) \cong \CHW^n(X,L).
		\]		
		
		\item For a morphism $f\colon Y\to X$ of smooth varieties over $k$ one has a functorial pullback homomorphism 
		\[
		f^*\colon \CHW^n(X,L) \to \CHW^n(Y,f^*L)
		\]
		given by the pullback in the cohomology of sheaves. Further, if $f$ is proper then one has a functorial pushforward homomorphism \cite[Section~2.3]{Fas20}
		\[
		f_*\colon \CHW_n(Y,f^*L) \to \CHW_n(X,L)
		\]
		induced by the transfers on Rost-Schmid complex. For a closed embedding $i\colon Z\to X$ of smooth varieties, with $j\colon X\setminus Z\to X$ the open embedding of the complement,   the localization sequence
		\[
		\CHW_n(Z,i^*L) \xrightarrow{i_*} \CHW_n(X,L) \xrightarrow{j^*} \CHW_n(X\setminus Z,j^*L)
		\]
		is exact \cite[Section~2.2]{Fas20}.		
		\item Let $F/k$ be a field extension of finite degree. Since $k$ is perfect, $F$ is separable over $k$, so    the field trace $\tr^F_k\colon F\to k$ is a nonzero $k$-linear functional. This gives rise to the Scharlau transfer 
		\[
		(\tr^F_k)_*\colon \GW(F)\to \GW(k),  
		\]
		an additive homomorphism, which is given on generators $\langle a \rangle\in\GW(F)$ by defining  $(\tr^F_k)_*(\langle a \rangle)$ to be the quadratic form $x\mapsto \tr^F_k(ax^2)$ on the $k$-vector space $F$. The pushforward homomorphism in Chow-Witt groups
		\[
		\pi_*\colon \CHW_0(\Spec F) \to \CHW_0(\Spec k)
		\]
		for the morphism $\pi\colon \Spec F\to \Spec k$ coincides by \cite[Example~1.23]{Fas20} with the Scharlau transfer $(\tr^F_k)_*$ under the identifications 
		\[
		\CHW_0(\Spec F)=K^{\MW}_0(F)\cong \GW(F),\quad \CHW_0(\Spec k)= K^{\MW}_0(k)\cong \GW(k).
		\]
		
		\item 
		For a rank $n$ vector bundle $E$ over a smooth variety $X$ over $k$ one has an \textit{Euler class} 
		\[
		e(E)=s^*s_*(1_X)\in \CHW^n(X,(\det E)^\vee)
		\]
		where $s\colon X\to E$ is the zero section. This class is natural with respect to pullbacks \cite[Proposition~3.1.1]{AF16}. 
		
		\item Let $E$ be a rank $n$ vector bundle over a smooth affine variety $X$ of dimension $n$ over $k$. Suppose that $\det E\cong \struct_X$. Then $e(E)=0$ if and only if $E$ has a non-vanishing section. For $n=1$ there is nothing to prove, for $n=2$ this was shown in \cite[Theorem~2.2]{BM00} and for general $n$ this follows from the results of \cite[Chapter~8]{Mor12}, \cite[Theorem~1]{AHW17} and \cite[Theorem~5.6]{AF16}, see Theorem~\ref{thm:Morel_criterion} below for the details.
		\end{enumerate}
\end{definition}

\begin{theorem}[Barge-Morel, Morel, Asok-Fasel, Asok-Hoyois-Wendt] \label{thm:Morel_criterion}
	Let $k$ be a perfect field and $E$ be a rank $n$ vector bundle over a smooth affine variety $X$ of dimension $n$ over $k$. Suppose that $\det E\cong \struct_X$. Then $E$ admits a non-vanishing section if and only if $e(E)=0$.
\end{theorem}
\begin{proof}
	In the case of $n=1$ there is nothing to prove, so assume $n\ge 2$. It follows from \cite[Theorem~8.14]{Mor12} that $E$ admits a non-vanishing section if and only if a certain obstruction-theoretic Euler class $e_\ob(E)\in \CHW^n(X)$ vanishes. Note that in {\it loc.cit.} it was assumed $n\ge 4$ because of the assumption $r\neq 2$ in \cite[Theorem~8.1.(3)]{Mor12}, which can be removed using \cite[Theorem~1]{AHW17}. It follows from \cite[Theorem~5.6]{AF16} that $e_\ob(E)=0$ if and only if $e(E)=0$, whence the claim.
\end{proof}

\begin{remark}
	We expect that Theorem~\ref{thm:Morel_criterion} generalizes to vector bundles with possibly non-trivial determinant and to general fields, removing the assumptions $\det E\cong \struct_X$ and $k$ being perfect.
\end{remark}

\begin{remark}\label{rem:StabTriv} Let $X$ be a smooth hypersurface in $\A^{n+1}_k$. Since $k[x_1,\ldots, x_{n+1}]$ is a UFD, the ideal of $X$ is principal, $I_X=(F)$. We have the tangent-normal sequence describing the tangent bundle $T_X$ as
\[
0\to T_X\to \struct_X^{n+1}\xrightarrow{\nabla F}\struct_X\to 0,
\]
where $\nabla F:=(\partial F/\partial x_1,\ldots,\partial F/\partial x_{n+1})$ is the usual gradient of $F$. Since $X$ is affine, this sequence splits, in particular, $T_X$ is stably trivial and $\det T_X\cong \struct_X$, so Theorem~\ref{thm:Morel_criterion} is applicable to $T_X$.
\end{remark}

\begin{definition} \label{def:odeg}
Let $X$ be a smooth proper variety over a perfect field $k$ with the structure morphism $\pi\colon X\to \Spec k$. Then the \textit{quadratic degree map}
\[
{\odeg}:=\pi_*\colon \CHW_0(X) \to \CHW_0(\Spec k) \cong \GW(k)
\]
is the pushforward homomorphism for the structure morphism. 
\end{definition}

\begin{lemma} \label{lem:birational_invariance}
	Let $Q$ be a smooth projective quadric over a perfect field $k$. Suppose that $Q(k)\neq \emptyset$. Then the quadratic degree map 
	\[
	{\odeg}\colon \CHW_0(Q) \to \GW(k)
	\]
	is an isomorphism.
\end{lemma}
\begin{proof}
	Let $n=\dim Q$ and consider the commutative square
	\[
	\xymatrix{
		Y \ar[d]_f \ar[r]^g & Q \ar[d]^\pi \\
		\PP^n \ar[r]^(0.4)p & \Spec k
	}	
	\]
	where $Q \xleftarrow{g} Y \xrightarrow{f} \PP^n$ is the resolution of the birational morphism $Q\dashrightarrow \PP^n$ given by the projection from a rational point on $Q$. Note that all maps in this square are proper.  This gives a commutative diagram of pushforward homomorphisms
	\[
	\xymatrix{
		\CHW_0(Y) \ar[d]_{f_*} \ar[r]^{g_*} & \CHW_0(Q) \ar[d]^{\pi_*}\\
		\CHW_0(\PP^n) \ar[r]^{p_*} & \GW(k)
	}	
	\]
	It follows from the birational invariance of $\CHW_0$ \cite[Corollary~2.2.11]{Fel22} that $f_*$ and $g_*$ are isomorphisms. Recall that $\omega_{\PP^n}=\struct_{\PP^n}(-n-1)$ whence
	\[
	\CHW_0(\PP^n)=\CHW^n(\PP^n,\struct(-n-1))\cong 
	\begin{cases}
		\CHW^n(\PP^n), & \text{$n$ is odd}, \\
		\CHW^n(\PP^n,\struct_{\PP^n}(-1)), & \text{$n$ is even}.
	\end{cases}
	\]
	The homomorphism $p_*$ is an isomorphism by \cite[Corollary~11.8]{Fas13}. Thus $\pi_* = \odeg$ is an isomorphism as well.
\end{proof}

\begin{definition}
A smooth projective scheme $X$ over $k$ has a {\em quadratic Euler characteristic}
$\chi(X/k)\in \GW(k)$, arising from the categorical Euler characteristic of the dualizable object $\Sigma^\infty_{\PP^1}X_+$ in the motivic stable homotopy category $\SH(k)$, together with Morel's theorem \cite[Theorem 6.4.1, Remark 6.4.2]{Mor04} identifying the endomorphisms of the unit in  $\SH(k)$ with $\GW(k)$ (see \cite[\S1]{Hoy14} and \cite[\S2.1]{Lev20} for details). The motivic Gau{\ss}-Bonnet theorem (\cite[Theorem 4.6.1]{DJK21}, \cite[Theorem 5.3]{LR20}) gives the identity 
\begin{equation}\label{eqn:GaussBonnet}
	\chi(X/k)=\odeg(e(T_X))\in \GW(k).
\end{equation}
	
\end{definition}

\begin{theorem} \label{thm:computational}
	Let $Q^o$ be the affine quadric over a perfect field $k$ given by the equation 
	\[
	a_1x_1^2+a_2x_2^2+\dots+a_{n+1}x_{n+1}^2=1
	\]
	with $a_1,\dots,a_{n+1}\in k^\times$ and let $Q^\infty$ be the projective quadric given by the equation
	\[
	a_1x_1^2+a_2x_2^2+\dots+a_{n+1}x_{n+1}^2=0.
	\]
	Then the following holds.
	\begin{enumerate}
		\item 
		If  $n$ is odd, then the tangent bundle $T_{Q^o}$ has a non-vanishing section.
		\item 
		If $n>0$ is even and the tangent bundle $T_{Q^o}$ has a non-vanishing section, then 
		\[
		\frac{n}{2}\cdot \langle 1, -1 \rangle + \langle 2, 2\cdot \prod_{i=1}^{n+1} a_i \rangle \in \odeg (\CHW_0(Q^\infty))\subseteq\GW(k).
		\]
		\item 
		If $n>0$ is even and $Q^o$ has a rational point, then the tangent bundle $T_{Q^o}$ has a non-vanishing section if and only if 		\[
		\frac{n}{2}\cdot \langle 1, -1 \rangle + \langle 2, 2\cdot \prod_{i=1}^{n+1} a_i \rangle \in \odeg (\CHW_0(Q^\infty))\subseteq\GW(k).
		\]
	\end{enumerate}
\end{theorem}

\begin{proof}
	(1)  We have settled the case of odd $n$ in Remark~\ref{rem:Explicit}.

	(2,3) Let $Q$ be the compactification of $Q^o$ given by the equation 
	\[
	a_1x_1^2+a_2x_2^2+\dots+a_{n+1}x_{n+1}^2=x_0^2
	\]
	in the projective space $\PP^{n+1}$ and let $j\colon Q^o\to Q$ be the open embedding. Then $Q^\infty=Q\setminus j(Q^o)$; let $i\colon Q^\infty\to Q$ be the closed embedding. Consider the localization sequence and the quadratic degree homomorphisms
	\[
	\xymatrix{
		\CHW_0(Q^\infty) \ar[r]^{i_*} \ar[dr]_{\odeg} & \CHW_0(Q) \ar[d]^{\odeg} \ar[r]^{j^*} & \CHW_0(Q^o) \\
		& \GW(k) &
	}
	\]
	
	We have the identifications
	 \[
	\CHW_0(Q)=\CHW^n(Q,\omega_Q)=\CHW^n(Q,(\det T_Q)^\vee), 
	\]
	so we consider the Euler class $e(T_Q)\in \CHW^n(Q,(\det T_Q)^\vee)$ as being in $\CHW_0(Q)$.
	
	Exactness of the localization sequence yields that the Euler class $e(T_{Q^0})= e(j^*T_{Q}) = j^* e(T_{Q})$ vanishes if and only if $e(T_{Q}) \in i_*\CHW_0(Q^\infty)$. By \eqref{eqn:GaussBonnet} and \cite[Corollary~12.2]{Lev20} we have
	\begin{equation}\label{eqn:EulerChar}
		\odeg(e(T_Q))= \chi(Q/k) =
		\frac{n}{2}\cdot \langle 1, -1 \rangle + \langle 2, 2\cdot \prod_{i=1}^{n+1} a_i \rangle.
	\end{equation}

	Suppose that $T_{Q^o}$ has a non-vanishing section. Then $e(T_{Q^0})=0$ and hence $e(T_{Q})$ is in $i_*\CHW_0(Q^\infty)$. Taking quadratic degrees and using formula~\ref{eqn:EulerChar} we obtain
	\[
	\odeg (e(T_{Q})) = \frac{n}{2}\cdot \langle 1, -1 \rangle + \langle 2, 2\cdot \prod_{i=1}^{n+1} a_i \rangle \in \odeg (i_*\CHW_0(Q^\infty)) = \odeg  (\CHW_0(Q^\infty)),
	\]
	proving (2) and one implication in (3).
	
	Now suppose that $Q^o$ has a rational point and $\frac{n}{2}\cdot \langle 1, -1 \rangle + \langle 2, 2\cdot \prod_{i=1}^{n+1} a_i \rangle \in \odeg  (\CHW_0(Q^\infty))$. Note that $\odeg(\CHW_0(Q^\infty)) = \odeg (i_*\CHW_0(Q^\infty))$, whence Lemma~\ref{lem:birational_invariance} combined with \eqref{eqn:EulerChar} show that $e(T_{Q})\in i_*\CHW_0(Q^\infty)$,  yielding $e(T_{Q^o})=0$. Remark~\ref{rem:StabTriv} provides an isomorphism $\det T_{Q^o}\cong \struct_{Q^o}$, whence Theorem~\ref{thm:Morel_criterion} implies that $T_{Q^o}$ has a non-vanishing section, completing the proof of (3).
\end{proof}

\begin{example}\label{exs:Examples_computational}
	Let $S^2_k$ be the quadric over a field $k$ given by the equation $x^2+y^2+z^2=1$.
	\begin{enumerate}
		\item If the equation $x^2+y^2=-1$ has a solution over $k$ then the conic $C_k\subseteq \PP^2_k$ given by the equation $x^2+y^2+z^2=0$ has a rational point whence $\odeg (\CHW_0(C_k))=\GW(k)$ and Theorem~\ref{thm:computational} yields that $T_{S^2_k}$ has a non-vanishing section. In particular \cite[Example~XI.2.4.(2) and~(6)]{Lam05} yield that this holds for $k=\QQ_p$ the field of $p$-adic numbers and for $k=\F_{p^n}$ a finite field,  with $p\neq 2$ in both cases. By base-change, the same follows if $k\supset \QQ_p$ or $k$ has characteristic $p$, with  $p\neq 2$. An explicit non-vanishing section of $T_{S^2_k}$ in these cases may be found as in Remark~\ref{rem:Explicit}.
		\item Let $k=\QQ_2$ be the field of dyadic numbers, then the equation $x^2+y^2=-1$ has no solution over $k$ by e.g. \cite[Example~XI.2.4.(7)]{Lam05} and the conic $C_k\subseteq\PP^2_k$ given by the equation $x^2+y^2+z^2=0$ has no rational points. However, it is clear that $C_k$ has a rational point over $\QQ_2(\sqrt{-2})$. Moreover, since $2$ is equivalent to $-14$ modulo squares in $\QQ_2$ (see, e.g., \cite[Corollary~VI.2.24]{Lam05}), 	
		$C_{k}$ has the point $(6+\sqrt{-14}, 10, -2+3\sqrt{-14})$ over $\QQ_2(\sqrt{2})$. A straightforward computation shows that $(\tr^{\QQ_2(\sqrt{\pm 2})}_{\QQ_2})_*(\langle 1\rangle)=\langle 2, \pm 1\rangle$ whence
		\[
		\langle 1, -1\rangle + \langle 2, 2\rangle = (\tr^{\mathbb{Q}_2(\sqrt{2})}_{\QQ_2})_*(\langle 1\rangle) + (\tr^{\QQ_2(\sqrt{-2})}_{\QQ_2})_*(\langle 1\rangle) \in \odeg(\CHW_0(C_k))
		\]
		and Theorem~\ref{thm:computational}.(3) yields that $T_{S^2_{k}}$ has a non-vanishing section. Alternatively, we have $(\tr^{\QQ_2(\sqrt{2})}_{\QQ_2})_*(\langle 1+\sqrt{2}\rangle)=\langle -2, 1\rangle$ whence 
\[
\langle 1, -1\rangle + \langle 2, 2\rangle = \langle 2, -2\rangle + \langle 1, 1\rangle =
(\tr^{\mathbb{Q}_2(\sqrt{2})}_{\QQ_2})_*(\langle 1, 1+\sqrt{2}\rangle)   \in \odeg(\CHW_0(C_k)).
\]
		Note that  $C_k$ does not have a $\QQ_2$-point, so we cannot apply Remark~\ref{rem:Explicit}. We  were not able to produce an explicit non-vanishing section in this case.
	\end{enumerate}
\end{example}

\begin{example}\label{exs:Examples_computational2} We can expand on the last example as follows.  Let $S^{2n}_k$ be the quadric over a field $k$ given by the equation $\sum_{i=1}^{2n+1}x_i^2=1$, $n>0$, and suppose $k$ contains a $p$-adic field $\mathbb{Q}_p$. Then $T_{S^{2n}_k}$ has a non-vanishing section. Indeed, letting $C^{2n}_k\subseteq\mathbb{P}^{2n}_k$ be the projective quadric defined by $\sum_{i=1}^{2n+1}x_i^2=0$, we have just seen that   $\langle 1, -1 \rangle + \langle 2, 2 \rangle$ is in  $ \odeg (\CHW_0(C^2_k))\subseteq \GW(k)$. But for an arbitrary quadratic extension $k(\sqrt{a})$ of $k$, we have
 \[
 (\tr^{k(\sqrt{a})}_k)_*(\langle\sqrt{a}\rangle)=\langle 1, -1 \rangle
 \]
 so $\langle 1, -1 \rangle$ is in $ \odeg (\CHW_0(Q))$ for {\em every} smooth projective
 quadric $Q$ over $k$, and thus we have
 \[
 n\cdot \langle 1, -1 \rangle + \langle 2, 2 \rangle \in \odeg (\CHW_0(C^2_k))\subseteq
  \odeg (\CHW_0(C^{2n}_k))\subseteq \GW(k).
  \]
  We then apply Theorem~ {\hyperref[thm:computational]{\ref*{thm:computational}.(3)}} to conclude that 
  $T_{S^{2n}_k}$ has a non-vanishing section. Just as in Example~ {\hyperref[exs:Examples_computational]{\ref*{exs:Examples_computational}.(1)}}, we can produce an explicit non-vanishing section if $C^{2n}_k$ has a $k$-rational point. This is the case if $\mathbb{Q}_p\subseteq k$ with $p$ an odd prime or if $n\ge 2$ and $\mathbb{Q}_2\subseteq k$ (for this last case, see \cite[Theorem~VI.2.12]{Lam05}).  
    \end{example}

\begin{remark}\label{rem:CHWS2Q2} Returning to the example of $S^2_{\QQ_2}$, Nanjun Yang asked if one could completely compute $\CHW_0(S^2_{\QQ_2})$. From the localization exact sequence in the proof of Theorem~\ref{thm:computational}, we have the isomorphism
\[
\CHW_0(S^2_{\QQ_2})\cong \GW(\QQ_2)/\odeg(\CHW_0(C_{\QQ_2})).
\]
We may identify $C_{\QQ_2}$ with the Severi-Brauer variety associated to the standard quaternions $\HH_{\QQ_2}$ over $\QQ_2$. By \cite[XIII, Proposition~6]{Serre}, the Brauer group of a local field $K$ is isomorphic to $\QQ/\Z$ by the map
\[
\text{inv}_K: \text{Br}(K)\to \QQ/\Z,
\]
so for a degree 2 extension $k\supset \QQ_2$, $C_{\QQ_2}(k)\neq\0$ if and only if $k$ splits $\HH_{\QQ_2}$, i.e., if and only if the invariant $\text{inv}_k(\HH_k)$ in $\QQ/\Z$ is zero. But by \cite[XIII, Proposition~7]{Serre}, 
\[
\text{inv}_k(\HH_k)=2\cdot \text{inv}_{\QQ_2}(\HH_{\QQ_2})=2\cdot \frac{1}{2}=0, 
\]
so $C_{\QQ_2}(\QQ_2(\sqrt{a}))\neq\0$  for every non-square $a\in \QQ_2^\times$. 

Thus, $\langle 2, 2a\rangle =\tr_{\QQ_2(\sqrt{a})/\QQ_2}(\langle 1\rangle)$ is in $\odeg(\CHW_0(C_{\QQ_2}))$ for all non-squares $a$, so the ideal $I$ in $\GW(\QQ_2)$ generated by the forms 
\[
\{\langle 2, u\rangle \mid u=3,5,7\}\cup \{\langle 2, 2u\rangle \mid u=1, 3,5,7\}.
\]
is contained in $\odeg(\CHW_0(C_{\QQ_2}))$.
It is easy to see that $I$ is exactly the ideal in $\GI(\QQ_2)$ of even rank forms in $\GW(\QQ_2)$; since $\odeg(\CHW_0(C_{\QQ_2}))$ is clearly contained $\GI(\QQ_2)$, we have  
$\odeg(\CHW_0(C_{\QQ_2}))=\GI(\QQ_2)$ and
\[
\CHW_0(S^2_{\QQ_2})\cong \Z/2
\]
via the mod 2 rank map $\GW(\QQ_2)\to \Z/2$.
\end{remark}

\section{Scharlau's transfer for closed points on a quadric}

\begin{definition}
	Let $F$ be a field, $\chark F\neq 2$. We denote $\GI(F)\subseteq \GW(F)$ the ideal consisting of the even dimensional (virtual) regular quadratic forms in the Grothendieck-Witt ring of $F$.
\end{definition}

\begin{definition} \label{def:Scharlau}
	Let $E/F$ be a field extension of finite degree, $\chark F\neq 2$, and $s\colon E\to F$ be a nonzero $F$-linear functional. Then the \textit{Scharlau transfer} \cite{Sch69}
	\[
	s_*\colon \GW(E) \to \GW(F)
	\]
	is the additive homomorphism such that $s_*(\langle a \rangle)$ is the quadratic form $x\mapsto s (ax^2)$ on the $F$-vector space $E$. See \cite[Chapter~VII]{Lam05} for some of the properties of the Scharlau transfer. 
	
	Let $X$ be a variety over $k$. The \textit{transfer ideal of $X$} is given by
	\[
	\GI^{\tr}_X= \sum_{x\in X_{(0)}} (s_x)_* (\GW(k(x))) \subseteq \GW(k)
	\]
	with the sum taken over all the closed points of $X$ and $\{s_x\colon k(x)\to k\}_{x\in X_{(0)}}$ being a chosen set of nonzero $k$-linear functionals. Note that $\GI^{\tr}_X$ does not depend on the choices of $s_x$ \cite[Remark~VII.1.6(C)]{Lam05}. It is easy to see that the transfer ideal admits an alternative description as 
	\[
	\GI^{\tr}_X= \sum_{\substack{\text{$F/k$ finite}\\ X(F)\neq \emptyset}} (s_F)_* (\GW(F)) \subseteq \GW(k)
	\]
	with the sum taken over all the (isomorphism classes) of field extensions $F/k$ of finite degree such that $X_F$ has a rational point and $\{s_F\colon F\to k\}_F$ being some chosen set of nonzero $k$-linear functionals.
\end{definition}

\begin{remark}
	The transfer ideal $s_* (\GW(F))\subseteq \GW(k)$ for an extension of fields $F/k$ of finite degree is a classical object of study, see e.g. \cite[Chapter~VII]{Lam05}. This agrees with the notion introduced above if one considers $F$ as a zero dimensional variety $\Spec F$ over $k$.
\end{remark}

\begin{lemma} \label{lem:transfer_ideal_as_degree}
	Let $X$ be a smooth proper variety over a perfect field $k$. Then
	\[
	\GI^{\tr}_X = \odeg(\CHW_0(X)).
	\]
\end{lemma}
\begin{proof} This follows from the description of $\CHW_0(X)$ via the cohomology of the Rost-Schmid complex and the fact that the pushforward in $\CHW_0$ for a separable field extension of finite degree coincides with the Scharlau transfer for the field trace.
\end{proof}

\begin{lemma} \label{lem:GW_GI}
	Let $E/F$ be a field extension of even degree, $\chark F\neq 2$, and $s\colon E\to F$ be an $F$-linear nonzero functional. Then $s_*(\GW(E))\subseteq s_*(\GI(E))+\langle 1, -1\rangle\cdot \GW(F)$.
\end{lemma}
\begin{proof}
	Note that the claim does not depend on the choice of $s$ (cf. \cite[Remark~VII.1.6(C)]{Lam05}). Without loss of generality we may assume that $E=F(\alpha)/F$ is a simple extension, then there is a nonzero functional $s$ such that $s_*(\langle \alpha\rangle)=\frac{[E:F]}{2}\cdot \langle 1, -1\rangle$ \cite[Theorem~VII.2.3]{Lam05}. The claim follows, since $s_*(\phi)=s_*(\phi+\langle \alpha\rangle)-\frac{[E:F]}{2}\cdot \langle 1, -1\rangle$.
\end{proof}

\begin{definition}
	Let $q$ be a regular quadratic form over $k$. Then we use the following notation:
	\begin{itemize}
		\item 
		$D(q)$ is the set of nonzero values of $q$,
		\item
		$D(q)^2=\{a\cdot b\mid a,b\in D(q)\}$ is the set of products of pairs of nonzero values of $q$,
		\item
		$[D(q)]$ and $[D(q)^2]$ are the multiplicative subgroups of $k^\times$ generated by the respective sets. 
	\end{itemize}
	Note that if $1\in D(q)$ then $[D(q)]=[D(q)^2]$. Since all the  sets introduced above are stable under multiplication by squares $(k^\times)^2$ we will sometimes abuse the notation and denote in the same way the corresponding subsets of $k^\times/(k^\times)^2$.
\end{definition}

\begin{definition}
	For a regular quadratic form $q$ over $k$ of dimension $n$ the \textit{signed discriminant} $\disc_\pm(q)$ is given by the formula 
	\[
	\disc_\pm(q):=(-1)^{\frac{n(n -1)}{2}}\det A_q,
	\]
	where $A_q$ is a symmetric matrix representing $q$. This gives a well-defined map
	\[
	\disc_\pm\colon \GW(k) \to k^\times/(k^\times)^2.
	\]
	When restricted to the ideal $\GI(k)$ this map becomes a homomorphism.
\end{definition}

\begin{lemma} \label{lem:Scharlau_quadric_disc}
	Let $Q$ be a smooth projective quadric over $k$ defined by a quadratic form $q$ and take $\phi \in \GI_Q^{\tr}$. Then $\disc_\pm(\phi)$ is in $[D(q)^2]$.
\end{lemma}
\begin{proof}
	We may assume $Q(k)=\emptyset$, otherwise $q$ is isotropic and $D(q)=k^\times$ whence there is nothing to prove. Springer's theorem \cite[Theorem~VII.2.7]{Lam05} yields that for every closed point $x\in Q_{(0)}$ the degree $[k(x):k]$ is even. Hence $\GI_Q^{\tr}\subseteq \GI(k)$, whence $\disc_\pm$ restricted to $\GI_Q^{\tr}$ is a homomorphism. Thus it is sufficient to check the claim for $\phi = s_*(\psi)$ with $\psi \in \GW(k(x))$ for a closed point $x\in Q_{(0)}$ and $s$ a chosen $k$-linear functional $s\colon k(x)\to k$. Furthermore, by Lemma~\ref{lem:GW_GI},  we may assume $\psi\in \GI(k(x))$. By \cite[Chapter~II, Theorem~5.12]{Sch85} we have
	\[
	\disc_\pm (s_*(\psi)) = N_{k(x)/k}(\disc_{\pm} (\psi)) \in k^\times/(k^\times)^2.
	\]
	The quadric $Q_{k(x)}$ has a rational point whence $q_{k(x)}$ is isotropic and $D(q_{k(x)})=k(x)^\times$, in particular, $\disc_{\pm} (\psi) \in D(q_{k(x)})$. Then Knebusch's norm principle \cite[Theorem~VII.5.1]{Lam05} implies $N_{k(x)/k}(\disc_{\pm} (\psi))\in [D(q)^2]$.
\end{proof}
\begin{remark} If $q$ is a Pfister form then the last result was obtained in \cite[Lemma 3.6]{BFS14}.
\end{remark}

\begin{lemma} \label{lem:Scharlau_quadric}
	Let $Q$ be a smooth projective quadric over $k$ defined by a quadratic form $q$. Then $\langle a, b\rangle \in \GI_Q^{\tr}$ if and only if $-ab\in [D(q)^2]$.
\end{lemma}
\begin{proof}
	Since  $\delta_\pm(\langle a, b\rangle )=-ab$, one implication follows from Lemma~\ref{lem:Scharlau_quadric_disc}. For the other implication, first note that we may assume $Q(k)=\emptyset$, since otherwise $\GI_Q^{\tr}= \GW(k)$ and there is nothing to prove. Then there exists a closed point $x\in Q_{(0)}$ such that $[k(x):k]=2$ and we may choose $\alpha\in k$ such that $k(x)\cong k(\sqrt{\alpha})$. Then for the $k$-linear functional 
	\[
	s\colon k(\sqrt{\alpha}) \to k,\quad s(1)=0,\, s(\sqrt{\alpha})=1,
	\]
	one has $s_*(\langle 1\rangle)= \langle 1,-1\rangle$ whence 
	\[
	\langle 1,-1\rangle\in \GI_Q^{\tr}.
	\]
Taking $c_1, c_2\in k^\times$, we have 
\[
	\langle 1,-c_1c_2 \rangle  = \langle c_1\rangle \langle 1, -c_2\rangle + \langle 1, -c_1\rangle - \langle 1, -1\rangle.
\]
Recalling that $\GI_Q^{\tr}$  is an ideal in $\GW(k)$, it follows that
\begin{equation}\label{eqn:ProductRule}
\hbox{$\langle 1, -c_1\rangle,  \langle 1, -c_2\rangle \in \GI_Q^{\tr}\Rightarrow 
 \langle 1, -c_1c_2\rangle \in \GI_Q^{\tr}$.
 }
\end{equation}
We claim that 
\begin{equation}\label{eqn:Claim}
c,d\in D(q)\Rightarrow \langle 1, -cd\rangle \in \GI_Q^{\tr}.
\end{equation}
Accepting our claim for the moment, write $-ab\in [D(q)^2]$ as a product, $-ab:=\prod_i a_ib_i$, with $a_i, b_i\in D(q)$. By \eqref{eqn:Claim}, we have $\langle 1, -a_ib_i\rangle\in \GI_Q^{\tr}$ for each $i$.  By \eqref{eqn:ProductRule}, it follows that $\langle 1, -\prod_ia_ib_i\rangle=\langle 1, ab\rangle$ is in 
$\GI_Q^{\tr}$. We proceed to prove \eqref{eqn:Claim}. 
	
	First suppose that $\dim Q =0$. Then we may assume $q=x_1^2-\alpha x_2^2$ and $Q\cong \Spec k(\sqrt{\alpha})$. A straightforward computation with the same functional $s$ as above shows that
	\[
	s_*(\langle w_1 + \sqrt{\alpha} w_2\rangle) = \langle 1, -(u_1^2-\alpha u_2^2)(v_1^2-\alpha v_2^2)\rangle
	\]
	for $w_1=(u_1v_1+\alpha u_2v_2)/(u_1v_2+u_2v_1)$ and $w_2=1$, which proves \eqref{eqn:Claim} in this case.
	
	Now suppose that $\dim Q\ge 1$. Let $q(u)=a$, $q(v)=b$ and choose some $w\neq 0$ such that $\psi_q(u,w)=\psi_q(v,w)=0$, where $\psi_q$ is the symmetric bilinear form associated to $q$. Put $c=q(w)$ and let $\alpha=-a/c$. Then 
	\[
	q(u+\sqrt{\alpha} w) = q(u) + \alpha q(w) + 2\sqrt{\alpha} \psi_q(u,w) = q(u) + \alpha q(w) =0.
	\]
	Thus there is a closed point $x\in Q_{(0)}$ such that $k(x)\cong k(\sqrt{-a/c})$. Then for the same functional $s$ as above one has
	\[
	s_*(\langle\sqrt{-a/c}\rangle) = \langle 1, -a/c\rangle \in \GI^{\tr}_Q.
	\]
	The same argument shows that $\langle 1, -b/c\rangle$ is in $\GI^{\tr}_Q$. Using \eqref{eqn:ProductRule}, we see that $\langle 1,-ab \rangle$   also belongs to $\GI^{\tr}_Q$ and \eqref{eqn:Claim} follows.
\end{proof}

\begin{remark} \label{rem:norms}
	Let $Q$ be a smooth projective quadric over a field $k$ defined by a quadratic form $q$. Then the group $[D(q)^2]$ coincides with the group of norms $N_Q(k)$ of $Q$, i.e. with the multiplicative subgroup of $k^\times$ generated by the norms $N_{F/k}(a)$ with $a\in F^\times$ and $F/k$ being an extension of fields of finite degree such that $Q_F$ has a rational point \cite[Lemma~2.2]{CTS93}.
\end{remark}

\begin{definition}
	For $a_1,a_2,\dots,a_n\in k^\times$ an \textit{$n$-fold Pfister form $\langle\langle a_1,a_2,\dots, a_n \rangle\rangle$} is the quadratic form $\prod_{i=1}^n \langle 1, -a_i\rangle$ of dimension $2^n$. A regular quadratic form $q$ over $k$ is called a \textit{Pfister neighbor} if there exists $a\in k^\times$ such that $\langle a\rangle \cdot q$ is a subform of an $n$-fold Pfister form with $2^{n-1}<\dim q$. Note that the Pfister form containing $\langle a\rangle \cdot q$ for a Pfister neighbor $q$ is unique \cite[Proposition~X.4.17]{Lam05}.
\end{definition}

\begin{lemma} \label{lem:Pfister_neighbor}
	Let $q$ be a Pfister neighbor over a field $k$ with the associated Pfister form $\phi$ and let $Q$ and $\Phi$ be the projective quadrics given by $q=0$ and $\phi=0$ respectively. Then $\GI_{Q}^{\tr} = \GI_{\Phi}^{\tr}$ and $[D(q)^2]=[D(\phi)^2]=D(\phi)$.
\end{lemma}
\begin{proof}
	Suppose $\phi$ is an $n$-fold Pfister form and $q$ has dimension $m>2^{n-1}$.  Let $a\in k^\times$ be such that $\langle a\rangle \cdot q$ is a subform of $\phi$. Since the quadrics associated to $q$ and $\langle a\rangle \cdot q$ are the same and $[D(q)^2]=[D(\langle a \rangle \cdot q)^2]$, we may assume that $q$ is a subform of $\phi$. We claim that for a field extension $F/k$, $Q(F)\neq \emptyset$ if and only if $\Phi(F)\neq \emptyset$. Indeed, since $q$ is a subform of $\phi$ then $Q\subseteq \Phi\subset \PP_k^{2^n-1}$, and, moreover, $Q=\Phi\cap L$, where $L$ is some codimension $2^n-m$ linear subspace of $\PP^{2^n-1}_k$. Thus if $Q(F)\neq \emptyset$ then $\Phi(F)\neq \emptyset$. Now let $F$ be such that $\Phi(F)\neq \emptyset$. Then $\phi_F$ is isotropic whence hyperbolic \cite[Theorem~X.1.7]{Lam05}, so $\Phi_F$ contains a linear subspace  $L'\subset \PP^{2^n-1}_F$ of dimension $2^{n-1}-1$. Letting $L_F\subset \PP^{2^n-1}_F$ be the base-extension of $L\subset \PP^{2^n-1}_k$ to $F$, we see that $Q_F$ contains the linear subspace $L'\cap L_F\subset \PP^{2^n-1}_F$ of dimension at least $2^{n-1}-1-(2^n-m)=m-2^{n-1}-1\ge0$, whence $Q(F)\neq \emptyset$.
 	
	The transfer ideals are generated by the Scharlau transfers for the field extensions $F/k$ of finite degree such $Q(F)\neq \emptyset$ (respectively, $\Phi(F)\neq \emptyset$), whence it follows from the above that $\GI_{Q}^{\tr} = \GI_{\Phi}^{\tr}$. Then the equality $[D(q)^2]=[D(\phi)^2]$ follows from Lemma~\ref{lem:Scharlau_quadric} because $[D(q)^2]$ and $[D(\phi)^2]$ coincide with the sets of signed discriminants of the binary forms from the respective transfer ideals. Since $1\in D(\phi)$ then $[D(\phi)^2]=[D(\phi)]$ and the last equality $[D(\phi)]=D(\phi)$ follows from \cite[Theorem~XI.1.1]{Lam05}.
	
	Alternatively, for the equality $[D(q)^2]=[D(\phi)^2]$ one could apply the description of these groups as the groups of norms \cite[Lemma~2.2]{CTS93}.
\end{proof}

\section{Non-vanishing vector fields on affine quadrics via groups of values}

Using the results of the previous section, we can reformulate Theorem~\ref{thm:computational} in a more manageable form.

\begin{theorem}\label{thm:hedgehog}
	Let $q=a_1x_1^2+a_2x_2^2+\dots+a_{n+1}x_{n+1}^2$ be a quadratic form over a perfect field $k$ with $a_1,\dots,a_{n+1}\in k^\times$, and let $Q^o$ be the affine quadric given by the equation
	\[
	a_1x_1^2+a_2x_2^2+\dots+a_{n+1}x_{n+1}^2=1.
	\]
	Then the following holds.
	\begin{enumerate}
		\item 
		If $n$ is odd then the tangent bundle $T_{Q^o}$ has a non-vanishing section.
		\item 
		If $n>0$ is even and the tangent bundle $T_{Q^o}$ has a non-vanishing section then 
		\[-\prod_{i=1}^{n+1} a_i \in [D(q)^2].\]
		\item 
		If $n>0$ is even and $Q^o(k)\neq \emptyset$ then the tangent bundle $T_{Q^o}$ has a non-vanishing section if and only if $-1 \in [D(q)]$.
	\end{enumerate}
\end{theorem}
\begin{proof} The case of odd $n$ follows from Theorem~{\hyperref[thm:computational]{\ref*{thm:computational}.(1)}}, so we assume $n>0$ is even. Let $Q^\infty\subseteq \PP^{n}$ be the quadric given by $q=0$. By Lemma~\ref{lem:transfer_ideal_as_degree} we have
	\[
	\frac{n}{2}\cdot \langle 1, -1 \rangle + \langle 2, 2\cdot\prod_ia_i \rangle \in \odeg(\CHW_0(Q^\infty)) \iff 
	\frac{n}{2}\cdot \langle 1, -1 \rangle + \langle 2, 2\cdot\prod_ia_i \rangle\in \GI_{Q^\infty}^{\tr}.
	\]
Note that one trivially has $1\in [D(q)^2]$. Thus, Lemma~\ref{lem:Scharlau_quadric} yields $\langle 1,-1\rangle \in \GI_{Q^\infty}^{\tr} $ and implies in addition that
\[
	\frac{n}{2}\cdot \langle 1, -1 \rangle + \langle 2, 2\cdot\prod_ia_i \rangle\in \GI_{Q^\infty}^{\tr} \iff  -\prod_{i=1}^{n+1}a_i \in [D(q)^2].
\]
Applying Theorem~{\hyperref[thm:computational]{\ref*{thm:computational}.(2)}} proves (2).

 For (3), note that $Q^o(k)\neq\emptyset$ if and only if $1\in D(q)$, which implies that $[D(q)]=[D(q)^2]$. Since each $a_i$ is in $D(q)$, we see that $-\prod_{i=1}^{n+1} a_i \in [D(q)^2]$ if and only if $-1\in [D(q)]$. Thus, $-1\in [D(q)]$ if and only if $\frac{n}{2}\langle 1, -1 \rangle + \langle 2, 2\cdot\prod_ia_i \rangle$  is in $\odeg(\CHW_0(Q^\infty))$, and then 
 Theorem~{\hyperref[thm:computational]{\ref*{thm:computational}.(3)}}  implies the claim. 
\end{proof}
\begin{definition}
	Let $F$ be a field. The \textit{level} of $F$, denoted $s(F)$, is the minimal integer $n$ such that $-1\in D(x_1^2+x_2^2+\dots+x_n^2)$. If no such $n$ exists then $s(F)=\infty$. The level of a field is either infinite or a power of $2$ \cite[Pfister's Level Theorem~XI.2.2]{Lam05}.
\end{definition}

\begin{corollary} \label{cor:hedgehog}
	Let $S^n_k$, $n\ge1$, be the affine quadric over a field $k$ given by the equation
	\[
	x_1^2+x_2^2+\dots+x_{n+1}^2=1.
	\]
	Then the tangent bundle $T_{S^n_k}$ has a non-vanishing section if and only if one of the following holds:
	\begin{enumerate}
		\item $n$ is odd, 
		\item $n>0$ is even and $s(k)\le 2n+1$.		
	\end{enumerate}
\end{corollary}
\begin{proof}
	First assume that $\chark k>2$. Then Remark~\ref{rem:Explicit} yields that the tangent bundle $T_{S^n_k}$ has a non-vanishing section for every $n\ge 1$. At the same time $-1$ is a square or a sum of two squares in a finite field, whence in $k$, thus $s(k)$ is $1$ or $2$ and, in particular, $s(k)\le 2n+1$ for $n\ge 1$. This yields the claim in the positive characteristic.
	
	Assume $\chark k=0$, in particular, $k$ is perfect. In view of Theorem~\ref{thm:hedgehog} we need to show that for even $n$ one has 
	\[
	-1\in [D(x_1^2+x_2^2+\dots+x_{n+1}^2)]
	\]
	if and only if $s(k)\le 2n+1$. Let $m$ be such that $2^{m-1}<n+1\le 2^m$. Then $x_1^2+x_2^2+\dots+x_{n+1}^2$ is a Pfister neighbor with the associated Pfister form $x_1^2+x_2^2+\dots+x_{2^m}^2$. Lemma~\ref{lem:Pfister_neighbor} yields
	\[
	[D(x_1^2+x_2^2+\dots+x_{n+1}^2)]=D(x_1^2+x_2^2+\dots+x_{2^m}^2).
	\]
	Thus $-1\in [D(x_1^2+x_2^2+\dots+x_{n+1}^2)]$ if and only if $s(k)\le 2^m$. The claim follows since $s(k)$ is a power of $2$.	
\end{proof}

\begin{example}\label{exs:Examples}
	Let $S^2_k$ be the quadric over a field $k$ given by the equation $x^2+y^2+z^2=1$. If $k=\mathbb{R}$ then $T_{S^2_k}$ has no non-vanishing sections since by a classical result of Poincar\'e the real vector bundle $T_{S^2_k}(\mathbb{R})$ has no non-vanishing continuous sections in the Euclidean topology  \cite[Theorem~6.5.5]{tD08}. More generally, Corollary~\ref{cor:hedgehog} yields that $T_{S^2_k}$ has no non-vanishing sections if and only if $s(k)\ge 8$ (including the case of $s(k)=\infty$). In particular, $T_{S^2_k}$ has a non-vanishing section for the following fields (cf. Example~\ref{exs:Examples_computational}):
	\begin{enumerate}
		\item $k$ a quadratically closed field,
		\item $k$ a field of characteristic $p>2$,
		\item $k$ a non-Archimedean local field,
		\item $k$ a purely imaginary number field.
	\end{enumerate}
	See \cite[Example~XI.2.4]{Lam05} for the relevant computations of $s(k)$.
\end{example}

\begin{corollary}\label{cor:CohDim2}
	Let $k$ be a perfect field of characteristic $\neq2$ such that every Pfister form of dimension $8$ is hyperbolic, i.e., that $I(k)^3=0$,\footnote{It follows from the Milnor conjecture \cite[Corollary 7.5]{Voev},  \cite[Theorem 4.1]{OVV},  \cite[Theorem 1.1]{RO}, that $I(k)^3=0$ is equivalent to $k$ being of $2$-cohomological dimension at most $2$.} and let $Q^o$ be the affine quadric over $k$ given by the equation
	\[
	a_1x_1^2+a_2x_2^2+\dots+a_{n+1}x_{n+1}^2=1
	\]
	with $a_i\in k^\times$. Suppose that $n$ is odd, or that $n>0$ is even and $Q^o(k)\neq \emptyset$. Then the tangent bundle $T_{Q^o}$ has a non-vanishing section.
\end{corollary}
\begin{proof}
	In view of Theorem~\ref{thm:hedgehog} it is sufficient to show that for an even $n>0$ one has $-1\in [D(\langle a_1,a_2,\dots,a_{n+1}\rangle)]$. We claim that already 
	\[
	-1\in [D(\langle a_1,a_2,a_3\rangle)^2]\subseteq [D(\langle a_1,a_2,\dots,a_{n+1}\rangle)^2] = 	 [D(\langle a_1,a_2,\dots,a_{n+1}\rangle)].
	\]
	Indeed, note that the quadratic form $\langle a_1,a_2,a_3\rangle$ is a Pfister neighbor since 
	\[
	\langle a_1 \rangle \cdot \langle a_1,a_2,a_3\rangle = \langle 1,a_1a_2,a_1a_3\rangle
	\]
	is a subform of the $2$-fold Pfister form $\langle\langle  -a_1a_2,-a_1a_3\rangle\rangle=\langle1,a_1a_2,a_1a_3,a_2a_3\rangle$. Then Lemma~\ref{lem:Pfister_neighbor} yields 
\begin{equation}\label{eqn:PfisterId}
	[D(\langle a_1,a_2,a_3\rangle)^2] = D(\langle 1,a_1a_2,a_1a_3,a_2a_3\rangle).
\end{equation}
	The form $\langle 1,a_1a_2,a_1a_3,a_2a_3,1\rangle$ is again a Pfister neighbor with the associated 3-fold Pfister form $\langle\langle -a_1a_2, -a_1a_3, -1\rangle\rangle$. The latter form is hyperbolic by the assumption whence its subform $\langle 1,a_1a_2,a_1a_3,a_2a_3,1\rangle$ is isotropic by the same dimension count argument as in the proof of Lemma~\ref{lem:Pfister_neighbor}. It follows that the equation
	\[
	x_1^2+a_1a_2 x_2^2+a_1a_3x_3^2+a_2a_3x_4^2+x^2_5=0
	\]
	has a solution over $k$. This means that
	\[
	-1\in D(\langle 1,a_1a_2,a_1a_3,a_2a_3\rangle)
	\]
	yielding by \eqref{eqn:PfisterId} that $-1\in [D(\langle a_1,a_2,a_3\rangle)^2]$ and the claim.
\end{proof}

\begin{corollary}\label{cor:NumberField} Let $k$ be a number field and let  $Q^o$ be the affine quadric over $k$ given by the equation
	\[
	a_1x_1^2+a_2x_2^2+\dots+a_{n+1}x_{n+1}^2=1
	\]
	with $a_i\in k^\times$. Suppose $Q^o(k)\neq \emptyset$.
	\begin{enumerate}
	\item If $n$ is odd, then $T_{Q^o}$ has a non-vanishing section.
	\item If $n>0$ is even, then $T_{Q^o}$ has a non-vanishing section if and only if for each real embedding $\sigma \colon k\to \mathbb{R}$, $\sigma(a_i)<0$ for some $i$.
	\end{enumerate}
\end{corollary}
\begin{proof}  The case of odd $n$ follows from Theorem~{\hyperref[thm:hedgehog]{\ref*{thm:hedgehog}.(1)}}. 

For even $n>0$, put $q=a_1x_1^2+a_2x_2^2+\dots+a_{n+1}x_{n+1}^2$. Let $\sigma \colon k\to \mathbb{R}$ be an embedding such that $\sigma(a_i)>0$ for all $i$ and let $q_{\mathbb{R}}$ be $q$ extended to $\mathbb{R}$ using this embedding. Then $[D(q_\mathbb{R})]\subseteq \mathbb{R}_{>0}$ whence $-1\not \in [D(q)]$. Thus Theorem~{\hyperref[thm:hedgehog]{\ref*{thm:hedgehog}.(3)}} yields one direction of the desired implication.

For the other direction it suffices to show that if $n> 0$ is even and for every embedding $\sigma \colon k\to \mathbb{R}$ one has $\sigma(a_i)<0$ for some $i$ then $-1\in [D(q)]$.  Note that the assumption that  $Q^o(k)\neq \emptyset$ implies that for each real embedding $\sigma$ of $k$,   there is a $j$ with $\sigma(a_j)>0$.

First assume $n\ge 4$. Let $v$ be a place of $k$ and consider the quadratic form $q+x_{n+2}^2=a_1x_1^2+a_2x_2^2+\dots+a_{n+1}x_{n+1}^2+x_{n+2}^2$ over $k_v$. If $v$ is a finite place, then $q+x_{n+2}^2$ is isotropic since every quadratic form of dimension $\ge 5$ is isotropic over a local field \cite[Theorem~VI.2.12]{Lam05}. If $v$ is an infinite place then the assumption $\sigma(a_i)<0$ for some $i$ implies that $q+x_{n+2}^2$ is isotropic. Then \cite[Hasse-Minkowski Principle~VI.3.1]{Lam05} implies that  $q+x_{n+2}^2$ is isotropic over $k$ whence $-1\in D(q)$.

Now assume $n=2$. The form $\langle a_1,a_2,a_3\rangle=a_1x_1^2+a_2x^2+a_3x^2$ is a Pfister neighbor with the associated Pfister form $\langle 1,a_1a_2,a_1a_3,a_2a_3\rangle$ and Lemma~\ref{lem:Pfister_neighbor} yields
\begin{equation}\label{eqn:PfisterNbr}
[D(\langle a_1,a_2,a_3\rangle)^2]=D(\langle 1,a_1a_2,a_1a_3,a_2a_3\rangle).
\end{equation}
Let $v$ be a place of $k$ and consider the form $q'=\langle 1,a_1a_2,a_1a_3,a_2a_3,1\rangle$ over $k_v$. As above, if $v$ is a finite place then the form $q'$ is isotropic by \cite[Theorem~VI.2.12]{Lam05}. If $v$ is a complex place then $q'$ is clearly isotropic. If $v$ is a real place with the real embedding $\sigma_v\colon k\to \mathbb{R}$ then as $Q^o(k)\neq\0$, we have $Q^o(k_v)\neq\0$, hence there is a $j$ such that $\sigma(a_j)>0$.  Combined with our assumption that $\sigma_v(a_i)<0$ for some $i$, we see that  at least one of $\sigma_v (a_1a_2)$, $\sigma_v (a_1a_3)$ and $\sigma_v (a_2a_3)$ is negative, whence $q'$ is isotropic over $k_v$. Then \cite[Hasse-Minkowski Principle~VI.3.1]{Lam05} implies that  $q'$ is isotropic over $k$ meaning that $-1\in D(\langle 1,a_1a_2,a_1a_3,a_2a_3\rangle)$. By \eqref{eqn:PfisterNbr},  we thus have
$-1\in [D(\langle a_1,a_2,a_3\rangle)^2]$ and the claim follows.
\end{proof}

\begin{corollary}\label{cor:LocalFieldEtc} Let $k$ be a field  of one of the following types:
\begin{enumerate}
\item a finite field $\F_{p^n}$, $p>2$,
\item  a non-Archimedean local  field of characteristic zero,
\item the perfection of a local field of characteristic $p>2$,
\item  the perfection of the function field of a curve over a finite field.
\end{enumerate}
Let  $Q^o$ be the affine quadric over $k$ given by the equation
	\[
	q:=a_1x_1^2+a_2x_2^2+\dots+a_{n+1}x_{n+1}^2=1
	\]
	with $a_i\in k^\times$. Suppose that $n>0$, and in case $n=2$ and $k$ is of type (2,3,4), suppose in addition that $Q^o$ has a $k$-rational point. Then $T_{Q^o}$ has a non-vanishing section.
\end{corollary}

\begin{proof} In all the above cases,  $k$ is a perfect field of cohomological dimension $\le 2$, and the result follows from Corollary~\ref{cor:CohDim2}, once we know that $Q^o(k)\neq\0$ in case $n\ge2$ is even. For $k$ of type (1) every regular quadratic form in $\ge3$ variables is isotropic \cite[Proposition~I.3.4]{Lam05} and for $k$ of type (2,3,4),  every regular quadratic form in $\ge5$ variables is isotropic \cite[Theorem~VI.2.12, Corollary~VI.3.5]{Lam05}; applying this to $q-x_0^2$  shows that $Q^o(k)\neq\0$ in all cases to be considered. \end{proof}

\end{document}